\documentclass[12pt,a4paper]{article}
\usepackage[utf8]{inputenc}
\pdfoutput=1
\usepackage[linesnumbered,lined,boxed,commentsnumbered]{algorithm2e}
\usepackage[text={6.5in,9.5in},centering]{geometry}

\usepackage{amsmath,amsfonts,amssymb, amsthm}
\usepackage{caption, subcaption}
\usepackage{graphicx, tasks, fancyvrb}
\usepackage{wrm_macros}
\usepackage{hyperref} 
\hypersetup{colorlinks=true, citecolor=blue, linkcolor=blue}

\author{J.A. Mackenzie\footnote{Corresponding author: Department of Mathematics and Statistics, University of Strathclyde, Livingstone Tower, 26 Richmond Street, Glasgow G1 1XH, Scotland} \and W.R. Mekwi\footnote{School of Computing, Engineering and Physical Sciences, University of the West of Scotland, Paisley Campus, Paisley PA1 2BE, Scotland}}
\title{An $hr$-Adaptive Method for the Cubic Nonlinear Schr\"{o}dinger Equation\footnote{© 2019. This manuscript version is made available under the CC-BY-NC-ND 4.0 license http://creativecommons.org/licenses/by-nc-nd/4.0/}}
\date{\today}
\begin{document}
\maketitle
\section*{Abstract}
The nonlinear Schr\"{o}dinger equation (NLSE) is one of the most important equations in quantum mechanics, and appears in a wide range of applications including optical fibre communications, plasma physics and biomolecule dynamics. It is a notoriously difficult problem to solve numerically as solutions have very steep temporal and spatial gradients. Adaptive moving mesh methods ($r$-adaptive) attempt to optimise the accuracy obtained using a fixed number of nodes by moving them to regions of steep solution features. This approach on its own is however limited if the solution becomes more or less difficult to resolve 
over the period of interest. Mesh refinement methods ($h$-adaptive), where the mesh is locally coarsened or refined, is an 
alternative adaptive strategy which is popular for time-independent problems. In this paper, we consider the effectiveness of a combined method ($hr$-adaptive) to solve the NLSE in one space dimension. 
Simulations are presented indicating excellent solution accuracy compared to other moving mesh 
approaches. The method is also shown to control the spatial error based on the user's input error tolerance. Evidence is also presented indicating second-order spatial convergence using a novel 
monitor function to generate the adaptive moving mesh.

\section*{Keywords}
Adaptivity, moving mesh method, $hr$-adaptivity, cubic nonlinear Schr\"{o}dinger equation

\section{Introduction}
The nonlinear Schr\"{o}dinger equation (NLSE) is an important model of mathematical physics and has numerous applications in the physical sciences such as in optical fibre communications, plasma physics, biological and atomic physics, biomolecule dynamics, hydrodynamics etc. (see e.g. \cite{taha,molbook} and the references therein). Its properties are also of great mathematical interest and solving it numerically has always proved to be challenging due to the solitons that appear (and disappear) in the solution. 

When written in a reference frame moving at the group velocity of the carrying wave, the cubic NLSE takes the form 
\beq \label{sch}
\left . \begin{array}{ll}
  i\psi_t + \psi_{xx} + q|\psi|^2\psi = 0, & \quad -\infty < x< \infty,\: 0 <  t\le T, \\
  \qquad \psi(x,0) = \psi_0(x), & \quad -\infty < x< \infty, 
\end{array} \right\} 
\eeq
where $\psi: \mathbb{R}\times\mathbb{R}^+ \rightarrow \mathbb{C}$, $i^2=-1$ and $q$ is a given real constant. This work only considers the focusing case for which $q$ is positive.

The analytical properties of (\ref{sch}) are well known (see \cite{hmm85, ssv:86}). In particular, the linear Schr\"{o}dinger equation
\beq \label{linsch}
i\phi_t + \phi_{xx} = 0, 
\eeq
provides a model for the propagation of dispersive waves given by 
\[
\phi(x,t) =  \exp[i(kx-W(k)t)]\:,
\]
where $W(k) = k^2$. The phase speed is defined by $W(k)/k$ and clearly depends on $k$. The wave is thus dispersive \cite{hmm85, ssv:86}. It can be shown that the solutions of (\ref{linsch}) have an amplitude which decays like $t^{-1/2}$ for $t,x \to \infty$ with $x/t$ fixed (see \cite[\S 11.3]{whitham}). 

The cubic term in (\ref{sch}) opposes dispersion and hence makes it possible for the NLSE to possess solutions where the competing forces of nonlinearity and dispersion balance each other exactly. It has been shown that, using the initial condition $\psi(x,0)=$ sech $x$, this is achieved when $q=2$ \cite{hmm85}. This balance happens for $q=8, 18,$ or in general when $q=2N^2$ for integer $N$. These states correspond to bound states of $N$ solitons or solitary waves. Solitons are formed when a certain balance between nonlinearity and dispersion is reached.

The pure initial-value problem (\ref{sch}) has been shown to possess an infinite number of conservation laws (or so-called invariants of motion) \cite{s84}. The most common of these are the charge $Q$, and the energy $E$ \cite{fpv95}, given by
\bea\label{h1}
Q &=& \int_{-\infty} ^\infty |\psi(x,t)|^2 \; dx\\ \label{h2}
E &=& \int_{-\infty} ^\infty \left( |\psi_x(x,t)|^2 -\frac{q}{2}|\psi(x,t)|^4\right) \; dx \;.
\eea
These conservation laws play an important role in the analysis and dynamics of the NLSE and many numerical schemes have been built \cite{del81, hmm85, sm83, taha, barletti2017} that attempt to conserve either or both of these quantities during simulations. Striving for methods to conserve energy is very important as conservation of energy implies the $L^2$-boundedness of the solution, thus preventing blow-up of the computed solution \cite{ssv:86}. However, although lack of exact (energy) conservation may lead to nonlinear blow-up, time integrators which are energy-conserving may still not perform desirably for this problem \cite{ssv:86}.

Adaptive mesh methods for the solution  of partial differential equations (PDEs) generally attempt to optimise the number and/or placement of mesh cells. These techniques have been proposed for a few decades now and can broadly be categorised as $h$-, $p$- and $r$- adaptive methods. For $p$-adaptive strategies, the order of polynomials which represent the solution locally, is varied in order to achieve better solution accuracy. Mesh refinement or $h$-adaptive methods are the most widely used adaptive methods. These methods start with a fixed number of mesh elements and then, based on some local measure of problem difficulty or a posteriori error estimate, elements are added or removed as necessary during the time integration of the problem. This approach is particularly useful when the solution is required to attain a specified accuracy. It is common to combine $h$- and $p$-adaptive strategies and this approach is well understood and is widely used. However, the method is complex and the underlying a posteriori estimates used to drive them can be rather difficult to obtain for strongly nonlinear problems \cite{budd09,mitchell2011}.

In contrast, moving mesh or $r$-adaptive methods generally aim to use a given number of mesh cells efficiently by prioritising areas of large solution variation and placing fewer mesh cells where the solution is smooth. This approach has been used with much success in a wide range of applications (see e.g. \cite{hawk:91, bmr:01, cao-2003, mek07, budd09} and references therein). When compared with a fixed uniform mesh approach, $r$-adaptive methods are able to maintain similar or achieve better levels of accuracy using fewer mesh points \cite{hawk:91}. The main shortcoming of an $r$-adaptive approach is that the solution accuracy is limited due to the fixed number of nodes. 

Moving mesh methods for time-dependent problems may be further divided into \emph{static} and \emph{dynamic} methods. With static methods the PDE is solved between time steps assuming that the mesh is fixed. The mesh is then modified and the physical solution is transferred between meshes using interpolation or an appropriate projection technique. For dynamic methods, the number and connectivity of mesh elements is kept fixed for most of the simulation and the mesh is moved at every time step. A mesh equation is usually used to compute node speeds in order to move the mesh. This has the advantage that one is not required to transfer the solution between meshes because the PDE is reformulated to take into account the fact that the nodes are moving \cite{zt02, tang-03, t04, tang-05, lee2015, gao2016}. Furthermore, these methods are thought to do a good job of reducing ``dispersive errors", a property that is useful for this problem. The experiments presented here employ a dynamic technique, using a moving mesh PDE (MMPDE) to drive the movement of the mesh. 

Even though numerous studies of the NLSE have been undertaken \cite{ssc:86, fpv95, bcr-99, cen:01, xie2009, wanhuang2015, barletti2017} these methods usually focus on the conservative properties of the NLSE or only apply $r$-refinement in order to track or resolve the emerging solitons. Some authors have employed a combined $hr$-refinement approach to solve problems from applications including combustion and fluid mechanics \cite{adje:86, li2001, ong13, hu2015}, elastostatics and fracture mechanics \cite{ammons-98},  heat transfer problems \cite{kita-2000} and ocean modelling \cite{piggott-2005}. In previous work \cite{mek07} we proposed an $hr$-algorithm and applied it to a number of one-dimensional time-dependent PDEs.  To our knowledge, few authors have employed an $hr$-refinement strategy for the NLSE. In combining these methods, we aim to efficiently track the travelling soliton(s) (mesh movement) while controlling solution accuracy appropriately using mesh refinement. 

The layout of this paper is as follows: section \ref{sec:num_meth} describes the components of the $hr$-adaptive method for the NLSE and presents the algorithm. Numerical experiments are reported and discussed in section \ref{sec:res}. Conclusions and future work are given in section \ref{sec:conc}.

\section{Numerical method}\label{sec:num_meth}
We employ an adaptive mesh technique for the NLSE which involves adjusting both the number and position of nodes of the mesh, a so-called $hr$-adaptive approach. In addition, we also vary the time-step size during the computation. This section gives details of how the algorithm has been implemented.

\subsection{Spatial discretisation}
The first step in approximating the solution is re-defining the pure initial value problem (\ref{sch}) as an initial-boundary value problem in $x_l\le x \le x_r$, since over the time interval $[0, T]$ under consideration, the solutions of (\ref{sch}) are  assumed to be negligibly small outside $[x_l, x_r]$. At these boundaries, it is customary to pose homogeneous Dirichlet or Neumann boundary conditions \cite{hmm85, rev:86, ssc:86}; our experiments use the former.

The complex function $\psi$ is first decomposed into its real and imaginary parts $u$ and $v$ respectively, resulting in the coupled system of PDEs
\beq \label{s1}
\left. \begin{array}{l}
u_t + v_{xx} + q(u^2+v^2)v = 0, \\
v_t - u_{xx} - q(u^2+v^2)u = 0,\end{array} \right\} \: x\in (x_l, x_r), \; t\in (0,T], 
\eeq
\beq \label{s2}
\left. \begin{array}{l}
u(x,0) = \psi_{0R}(x),\quad v(x,0) = \psi_{0I}(x), \quad x\in (x_l, x_r), \\
u(x_l,t)= u(x_r,t) =  0,   \\
v(x_l,t)= v(x_r,t) =  0, \quad  t\in (0, T],  \end{array} \right\}
\eeq
where $\psi_{0R}$ and $\psi_{0I}$ are, respectively, the real and imaginary parts of
$\psi_0$. 

To include node movement, (\ref{s1}) is rewritten with respect to a moving reference frame and the problem is recast in terms of the independent variables $\xi$ and $t$, using the transformation 
\beq\label{mesh1}
x = x(\xi,t), \quad \xi \in \Oma_c \equiv [0,1], \quad t\in (0,T],
\eeq
from computational space $\Oma_c\times (0,T]$ to physical space $\Oma_p \times
  (0,T]$. A uniform mesh covering $\Oma_c$ is given by
\beq\label{xi1}
\xi_i = \frac{i}{N}, \quad i=0,1,\ldots,N,
\eeq
and the corresponding (nonuniform) mesh on $\Oma_p$ is
\beq\label{pi1}
x_l=x_0 < x_1(t)< \dots < x_{N-1}(t) < x_N=x_r,
\eeq
where
\beq\label{mesh2}
x_i(t) = x(\xi_i,t), \quad i=0,1,\ldots,N.
\eeq
In what follows, we define $\vecb{x}(t)=\left\{ x_i(t)\right\}_{i=0}^N$.

To incorporate mesh movement, it is convenient to express the Eulerian time derivatives in (\ref{s1}) in terms of the moving reference frame. The resulting equations are given by
\bea\label{u_dot}
\dot{u} &=& \dot{x} u_x - v_{xx} - q(u^2+v^2)v, \\
\label{v_dot} 
\dot{v} &=& \dot{x} v_x + u_{xx} + q(u^2+v^2)u,
\eea
where $\dot{u}$, $\dot{v}$ and $\dot{x}$ denote derivatives with respect to $t$ where $\xi$ is held constant. 

After applying central differencing for the spatial derivatives, the semi-discrete system on a moving grid is given by
\bea
\dot{U}_i &=& \dot{x}_i\frac{U_{i+1}-U_{i-1}}{h_{i+1}+h_i} - \frac{2}{h_{i+1}+h_i}
\left(\frac{V_{i+1}-V_i}{h_{i+1}} - \frac{V_i-V_{i-1}}{h_i} \right) \nonumber\\
& & -\: q(U_i^2 + V_i ^2)V_i, \label{usch}\\
\dot{V}_i &=& \dot{x}_i\frac{V_{i+1}-V_{i-1}}{h_{i+1}+h_i} + \frac{2}{h_{i+1}+h_i}
\left(\frac{U_{i+1}-U_i}{h_{i+1}} - \frac{U_i-U_{i-1}}{h_i} \right) \nonumber\\
& & +\: q(U_i^2 + V_i ^2)U_i, \label{vsch}
\eea
for $i=1,\ldots,N$, where $h_i = x_i - x_{i-1}$, and $U_i(t)$, $V_i(t)$ are the approximations to  $u(x_i ,t)$ and $v(x_i,t)$. Further, we define $\mathbf{U}=(U_0, U_1,\ldots,U_{N-1}, U_{N})^T$ and $\mathbf{V}=(V_0, V_1,\ldots,V_{N-1}, V_{N})^T$.

Application of the zero Dirichlet boundary conditions yields
\[
U_0 = 0,\:\: U_N = 0,\:\: V_0 = 0 \quad \text{ and } \quad V_N = 0. 
\]
The discretisation of $\dot{x}$ is described in section \ref{sec:ats}.

\subsection{Mesh adaptation}
The mesh adaptation process is a twofold process involving node movement and mesh refinement. 
\subsubsection{Mesh movement}
The aim of node movement is to place more nodes where the solution has features which are difficult to resolve. We  achieve this using a monitor function. This is typically a function of the solution gradient or curvature, and its choice is crucial to the success of the mesh movement process \cite{huang}. The placing of nodes is based on an equidistribution principle, where nodes are located to ensure that some measure of error is equally distributed over the mesh \cite{boo73}. 

The mesh velocity $\dot{x}$ is computed using a mesh generating equation. We use a moving mesh PDE (MMPDE) to achieve this. The impact of MMPDEs on the moving mesh method has been studied in some detail by Huang and his collaborators (\cite{hrr94,cao:01,huan:99,hr01}). Our experiments here use a one-dimensional version of one of the two-dimensional MMPDEs proposed by Huang and Russell \cite{huan:99} given by:
\beq\label{mmpde}
\ptl{x}{t} = \frac{1}{\tau}\left( M\ptl{x}{\xi}\right) ^{-2} \ptl{{}}{\xi}\left(
M\ptl{x}{\xi}\right), \:\:\xi \in \Oma_c, \:\: x(0,t) = x_l, \; x(1,t) = x_r,
\eeq
where $\tau$ is a temporal smoothing parameter and $M({\nu}(x,t))$ is a positive monitor function. The initial condition $x(\xi,0)$ is obtained by equidistribution of the a monitor function based on the initial condition at $t=0$. 

Recognising that the NLSE is a system of equations for the real and imaginary components $u$ and 
$v$, and motivated by the work in \cite{bm00:1, ver:89}, our experiments utilise the second derivative monitor function 
\begin{equation}
M=\frac{1}{2}(M_{u}+M_{v}), 
\label{2dev}
\end{equation}
where 
\begin{equation}
M_{u}=\varphi_{u}(t)+\sqrt{u_{xx}(x,t)}, \quad\quad \varphi_{u} (t) = \frac{1}{|\Oma_p|}\int_{\Oma_p} \sqrt{u^{2}_{xx}(x,t)} \: {\rm d}x
\label{eq:mu}
\end{equation}
and 
\begin{equation}
M_{v}=\varphi_{v}(t)+\sqrt{v_{xx}(x,t)}, \quad\quad \varphi_{v} (t) = \frac{1}{|\Oma_p|}\int_{\Oma_p} \sqrt{v^{2}_{xx}(x,t)} \: {\rm d}x.
\label{eq:mv}
\end{equation}
The role of $\varphi _{u,v}(t)$ is to ensure that some mesh points are placed in regions that do not have steep solution features e.g. outside boundary and interior layers. Some analysis on the effect of using monitor functions in general has been carried out 
in \cite{bm00:1, bmrs:01,huan:94, hr01}. 

Discretisation of (\ref{mmpde}) requires the evaluation of the monitor function (\ref{2dev}). If we let
\beq
g_{i+\FR{1}{2}} = \frac{U_{i+1} - U_i}{h_{i+1}}, 
\eeq
then an approximation $w_{i}$ to $\sqrt{|u_{xx}|}$ is given by
\beq\label{nuxx}
w_{i}= \sqrt{2\left |\left( \frac{g_{i+\FR{1}{2}} -g_{i-\FR{1}{2}}}{h_{i+1} + h_{i}} \right)\right |}.
\eeq

The floor, $\varphi_{u}$, on the monitor function, $M_{u}$, is given by a quadrature approximation
\beq\label{a_approx}
\varphi_{u} = \frac{1}{x_R - x_L} \sum_{i=0} ^{N-1}\frac{h_{i+1}}{2} (w_{i+1} + w_{i}).
\eeq

The discrete approximation of the monitor function (\ref{eq:mu}) is then given by
\beq
(M_u)_{i+\FR{1}{2}} = \varphi_{u} + \frac{1}{2} (w_{i+1} + w_{i}),\:\quad\quad i=0,1,\ldots, N-1.
\eeq
The approximation of $(M_v)_{i+\FR{1}{2}}$ is calculated similarly and finally 
\[
M_{i+\FR{1}{2}}=\frac{1}{2}((M_u)_{i+\FR{1}{2}}+(M_v)_{i+\FR{1}{2}}).
\]
To improve the robustness of the moving mesh method a smoothed monitor function is obtained as in \cite{mqs97} by setting
\beq\label{smooth}
\tilde{M}_{i+\frac{1}{2}} = \frac{{\ds \sum_{k=i-p}^{i+p} M_{k+\frac{1}{2}}\left(\frac{\gma}{\gma+1}\right)^{|k-i|}}}
      {{\ds \sum_{k=i-p}^{i+p} \left(\frac{\gma}{\gma+1}\right)^{|k-i|}}}\:, 
\eeq
with $\gma=2$ and $p=3$ in all our experiments. In addition, the summations only use terms that are well defined $(0\leq k \leq N-1)$. Using second-order central differences for the spatial derivatives in (\ref{mmpde}), we obtain a semi-discrete system of moving mesh equations defined by
\beq\label{x_i}
\dot{x}_i = \frac{4}{\tau}(\tilde{M}_i(h_{i+1}+h_i))^{-2}
(\tilde{M}_{i+\frac{1}{2}}h_{i+1}-\tilde{M}_{i-\frac{1}{2}}h_i),
\eeq
for $i=1,2,\ldots,N-1$, with $x_0 = x_l$ and $x_N = x_r$. The term $\tilde{M}_i$ in (\ref{x_i}) is given by 
\beq\label{m}
\tilde{M}_i =
\frac{\tilde{M}_{i-\frac{1}{2}}(x_{i+\frac{1}{2}}-x_i)+
  \tilde{M}_{i+\frac{1}{2}}(x_i-x_{i-\frac{1}{2}})}{x_{i+\frac{1}{2}} -
  x_{i-\frac{1}{2}}}\,,
\eeq
where 
\[ x_{i+\frac{1}{2}} = \frac{x_{i+1} + x_i}{2}\,.\]

\subsubsection{Mesh refinement}
In choosing the number of nodes to use at each time step, the objective is to keep some measure of the spatial solution error within predetermined bounds. To implement this strategy, we define an error monitor, which broadly measures global difficulty of the problem under consideration. 
Let us define
\beq\label{eta}
\eta(t) = \left ( \frac{1}{N} \int_{\Oma_p} M \; dx \right) ^2\:,
\eeq
where $M$ is given by (\ref{2dev}). If RTOL is a user-prescribed tolerance and $\ala$ and $\bta$ are given such that $\ala>1$ and $0<\bta<1$, then we would like to ensure that
\beq\label{ab}
\bta\textrm{ RTOL} \le \eta(t) \le \ala\textrm{ RTOL}\;.
\eeq
Note that $\eta(t)$ is simply a cheap, heuristic means that enables us to determine when to change the number of nodes. Also, the number of remeshings performed during a calculation is inversely proportional to the difference between $\ala$ and $\bta$.

To calculate the number of nodes required at time $t^{n+1}$ we first evaluate 
\beq\label{n}
\widetilde{N}^{n+1} = N^n \times \min \left( \text{maxfac},\max
\left[ \text{minfac},\kap\left(\frac{\eta(t)}{RTOL}\right)
  ^{\FR{1}{2}} \right] \right)\;,
\eeq
where maxfac, minfac and $\kap$ are user-defined parameters. The exponent of $\frac{1}{2}$ in the formula is related to the assumption that the scheme is second-order accurate in space. In our experiments we have set these values to be: maxfac $=2.0$, minfac $=1.2$ or $0.3$  for mesh enrichment or mesh coarsening respectively, and $\kap \ge 1$. The new number of mesh nodes is then 
\beq \label{n2}
N^{n+1} := \lfloor\widetilde{N}^{n+1}\rfloor + 1\:,
\eeq
where $\lfloor \cdot \rfloor$ denotes ``integer part of''.
The prediction of the number of nodes (\ref{n}) has been motivated using the same idea as has been used to control the time step in (\ref{dtsol}).

\subsection{Temporal adaptation}\label{sec:ats}
Integration from $t=t^n$ to $t=t^{n+1}$ is accomplished following the approach of Beckett \ea \cite{bmrs:01}. The semi-discretised systems (\ref{usch}), (\ref{vsch}) and (\ref{x_i}) for the solution and the mesh are clearly coupled. For efficiency, we solve them using a decoupled iterative procedure which we outline here.

We discretise (\ref{x_i}) in time using the first-order backward Euler method. This equation  is the governing equation for $\vecb{x}\equiv \{x_i\}_{i=0}^N$, for which $\tilde{M}_{i\pm\frac{1}{2}}$ and $(\tilde{M}_i(h_{i+1}+h_i))^{-2}$ are known for $i=1,2,\ldots,N-1$. We use an iterative approach where, initially, $\tilde{M}_{i\pm\frac{1}{2}}$ and $(\tilde{M}_i(h_{i+1}+h_i))^{-2}$ are evaluated at $t=t^n$. Subsequent iterations use values from the preceding cycle. The linear system which arises at each iteration is solved directly. The final mesh that is used is obtained by combining this mesh and the mesh from the previous iteration using under-relaxation so that the $(\nu + 1)$th iterate is given by 
\beq\label{xnu}
\vecb{x}^{[n+1,\nu +1]} = (1-\oma)\vecb{x}^{[n+1,\nu+1]}_* + \oma\vecb{x}^{[n+1,\nu]}\; ,
\eeq
where $\vecb{x}^{[n+1,\nu +1]}_*$ is the mesh obtained from solving (\ref{x_i}) and $\oma = 0.8$ in all our computations.

When integrating from $t=t^n$ to $t=t^{n+1}$, the system (\ref{usch}) and (\ref{vsch}) is considered a system of equations for the approximation of $\mathbf{U}$ and $\mathbf{V}$, where the node locations are available at $t=t^n$ and $t=t^{n+1}$. The vector $\dot{\vecb{x}}$ is replaced by $(\vecb{x}^{n+1} - \vecb{x}^n)/\Dta t^n$ and $\vecb{x}$ is evaluated in $[t^n,t^{n+1}]$ using the linear interpolant 
\[
\vecb{x}(t):= \vecb{x}^n + \dot{\vecb{x}}(t-t^n), 
\]
where $\vecb{x}^n$ is the approximation to $\vecb{x}$ at $t=t^n$. We also note that, while some authors (see e.g. \cite{cen:01}) use an upwinding approach for $\dot{x}$ in (\ref{u_dot}) and (\ref{v_dot}), such treatment is not used here.

The equations \eqref{usch} and \eqref{vsch} can be written as a coupled set of ODEs 
\beq\label{vdf}
\vecb{\dot{w}} = \vecb{f}(t,\vecb{w}),
\eeq
where $\vecb{f}: \mathbb{R} \times \mathbb{R}^{2(N+1)} \to \mathbb{R}^{2(N+1)}$ and 
\[
\vecb{w}=\left( \begin{array}{c}\mathbf{U} \\ \mathbf{V} \end{array} \right).
\]

The solution of (\ref{vdf}) at $t=t^{n+1}$ is computed using a second-order singly diagonally implicit Runge-Kutta method (SDIRK2) \cite{bmrs:01} which possesses excellent stability properties. The Butcher array is shown in Table~\ref{butcher}, where $\gma = (2-\sqrt{2})/2$.
\begin{table}[!ht]
\centering
\begin{tabular}{r|ll} 
$\gma$ & $\gma$ & $0$ \\
$1$ & $1-\gma$ & $\gma$ \\ \hline
\, & $1-\gma$ & $\gma$ \\ 
\end{tabular}  
\caption{Butcher array for SDIRK2 method.}
\label{butcher}
\end{table}
Integration from $t=t^n$ to $t=t^{n+1}=t^n+\Dta t^n$ is given by the solution of
\bea\label{rk}
\vecb{k}_1 &=& \vecb{f}(t^n + \gma\Dta t^n,\vecb{w}^n+\gma\Dta t^n\vecb{k}_1), \nonumber \\
\vecb{k}_2 &=& \vecb{f}(t^n + \Dta t^n,\vecb{w}^n + (1-\gma)\Dta t^n\vecb{k}_1 + \gma\Dta
t^n\vecb{k}_2),\\
\vecb{w}^{n+1} &=& \vecb{w}^n + \Dta t^n ((1-\gma)\vecb{k}_1 + \gma\vecb{k}_2),\nonumber
\eea
where $\vecb{w}^n$ denotes the value of $\vecb{w}$ at $t=t^n$. 

The values of $\vecb{k}_1$ and $\vecb{k}_2$ are obtained using a Newton iteration, terminating when convergence is achieve in the $l_\infty$ norm to within a user-prescribed tolerance, KTOL. 

Our algorithm also employs adaptive time stepping. To implement this, we use an error indicator based on the computed solution and a first order solution that can be cheaply obtained as a by product of the SDIRK computation used. Indeed, if $\hat{\vecb{w}}^{n+1}$ is a first order approximation to $\vecb{w}$ at time $t=t^{n+1}$, we obtain $\hat{\vecb{w}}^{n+1}$ using
\beq
\hat{\vecb{w}}^{n+1} = \vecb{w}^n + \Dta t^n \vecb{k}_1,
\eeq
and the error indicator ERR is calculated using the mesh-dependent $L_2$ error measure
\beq\label{err2}
\text{ERR } = \left(\sum_{i=0}^{N-1} (x_{i+1}^{n+1} - x_i^{n+1}) \left(\frac{ ||\mathbf{e}_i^{n+1}|| +
  ||\mathbf{e}_{i+1}^{n+1}||}{2}\right)^2 \right)^{\frac{1}{2}}\;,
\eeq
where 
\beq\label{er}
\mathbf{e}_i^{n+1} = \left (\begin{array}{c} U_i^{n+1} - \hat{U}_i^{n+1} \\ 
V_{i}^{n+1}-\hat{V}_{i}^{n+1} \end{array} \right ) \:, \text{ for } i=0, \ldots, N.
\eeq
In the event that ${\rm ERR}>{\rm ETOL}$, where ETOL is a user-supplied tolerance, the time step is 
rejected and repeated with a smaller time step. On the other hand if ${\rm ERR}<{\rm ETOL}$, then the predicted suitable time step, based on the solution, for the next time step is given by the formula (see \cite{hw2}) 
\beq\label{dtsol}
\Dta t^{n+1}_{(sol)} = \Dta t^n\times \textnormal{min} \left(
\textnormal{maxfac,\;max} \left[ \textnormal{minfac, } \;\eta \left(
\frac{\textnormal{ETOL}}{\textnormal{ERR}}\right) ^{\frac{1}{2}} \right]
\right),
\eeq
where typically $1.5 \leq$ maxfac $\leq 3$, $\eta \sim 0.6,$ and minfac$=0.1$.

A similar computation based on the mesh is carried out to predict the time step. Here, we use an indicator of grid accuracy, mesherr, given by
\beq\label{mesherr}
\textnormal{mesherr} = \| \vecb{x}^{[n+1,\upsilon]} - \vecb{x}^{[n+1,\upsilon-1]} \|_{l_{\infty}}, 
\eeq
where $\vecb{x}^{[n+1,\upsilon]}$ is the final mesh at $t=t^{n+1}$ obtained from \eqref{xnu}, and $\vecb{x}^{[n+1,\upsilon-1]}$ denotes the mesh at the previous pass in the iterative procedure. If MESHTOL is a user-supplied mesh tolerance and 
${\rm mesherr}>{\rm MESHTOL}$, then the time step is rejected and repeated with a smaller time step. 
If ${\rm mesherr}<{\rm MESHTOL}$, then the predicted mesh time step is calculated using:
\beq\label{dtmesh}
\Dta t^{n+1}_{(mesh)} = \Dta t^n\times \textnormal{min} \left(
\textnormal{maxfac,\;max} \left[ \textnormal{minfac,}\; 
\frac{\log(\textnormal{mesherr})}{\log(\textnormal{MESHBAL})} \right]
\right),
\eeq
where MESHBAL is a user-chosen parameter and we require that MESHBAL $<$ MESHTOL. 
We finally take the minimum of the solution time step \eqref{dtsol} and the predicted mesh time step \eqref{dtmesh} so that 
\beq\label{dt}
\Dta t^{n+1}=\min\left(\Dta t^{n+1}_{(sol)},\;\Dta t^{n+1}_{(mesh)}\right)\;.
\eeq

\subsection{The complete $hr$-algorithm}
The algorithm which we have used to solve the NLSE is shown in Algorithm \ref{talgo} below, where $\psi^{n}$ represents the numerical solution at time $t^n$. Note that, in the absence of steps \ref{starthref}-\ref{endhref}  and the computation of $N^0$ in step \ref{initstep}, the method performs $r$-adaptation with the initial number of nodes provided. 

The initial mesh $N^0$ is computed automatically so that $\eta(0)$ satisfies (\ref{ab}). This is achieved by equidistribution of the initial condition using the de Boor algorithm \cite{boo73} which
iterates until the difference between successive meshes is less than some prescribed tolerance, GTOL. 
\begin{algorithm}
\DontPrintSemicolon
\BlankLine
Initialise variables; select all parameters and $\Dta t^0$\;
Determine $\vecb{x}^0$ and $N^0$ based on initial condition, $\psi_0$ e.g. \eqref{ic1}\; \nllabel{initstep}
Set $n:=0$ and $t^1:=\Delta t^{0}$\;
\While{$t^{n+1}\leq T$}{
Set $\nu :=0$\;\nllabel{muo}
Set $\psi^{[n+1,0]}:=\psi^{n}$\;
\While{$\nu<4$}{
Solve the MMPDE for $\vecb{x}^{[n+1,\nu+1]}$\quad\eqref{xnu}\;
Using the SDIRK2 scheme find $\psi^{[n+1,\nu+1]}$\quad\eqref{rk}\;
$\nu:=\nu+1$\;
}
Calculate ERR and mesherr\quad\eqref{err2}, \eqref{mesherr}\;
\eIf{${\rm mesherr}<{\rm MESHTOL}$ {\rm and} ${\rm ERR}<{\rm ETOL}$}
{Determine $\Dta t^{n+1}$\quad\eqref{dt}}\;
{$\Delta t^{n}:=\Delta t^{n}/2$\;
Goto \ref{muo}\;}
Calculate $\eta(t^{n+1})$ using $\vecb{x}^{[n+1,4]}$ and $\psi^{[n+1,4]}$\quad\eqref{eta}\;
\eIf{\nllabel{starthref}$\beta\;{\rm RTOL}<\eta(t^{n+1})<\alpha\;{\rm RTOL}$}{$N^{n+1}:=N^{n}$\;}
{Compute $N^{n+1}$ using (\ref{n}) and (\ref{n2})\;
Generate $\ds \vecb{x}^{n+1}_{N^{n+1}}$\;
Transfer the solution $\psi^{[n+1,4]}$ onto $\ds \vecb{x}^{n+1}_{N^{n+1}}$\;}\nllabel{endhref}
Set $n:=n+1$\;
}
\caption{The $hr$-algorithm.}
\label{talgo}
\end{algorithm}

In the event that $h$-refinement is performed, the new grid $\ds \vecb{x}^{n+1}_{N^{n+1}}$ 
is obtained via equidistribution of the monitor function defined on the mesh $\ds \vecb{x}^{n+1}_{N^n}$. 
Using cubic interpolation, the solution is then transferred from $\ds \vecb{x}^{n+1}_{N^{n}}$ onto the newly computed grid, $\vecb{x}_{N^{n+1}}^{n+1}$. 
\section{Numerical experiments}\label{sec:res}
We now consider the performance of the $hr$-algorithm on three increasingly difficult test problems. 
\begin{table}[!ht]
\centering
{\renewcommand{\arraystretch}{1.3}
\begin{tabular}{|l|l|}
\hline
NHR & Number of $h$-refinements performed \\
NMAX & Largest number of nodes used \\
NMIN & Least number of nodes used \\
NSTP & Number of time steps used \\
JACS & Number of Jacobian computations throughout the simulation \\
& (when computing $\vecb{k}_1$ and $\vecb{k}_2$)\\
BS & Number of back solves needed when an $LU$ \\
& factorisation is used in the quasi-Newton iterations.\\
ETF & Number of times the time step is halved due to fail in error test i.e. \\
& ERR $>$ ETOL or mesherr $>$ MESHTOL (\S \ref{sec:ats}) \\
CTF & Number of convergence test fails in Newton iteration i.e. \\
& $\vecb{k}_1$ or $\vecb{k}_2$ fail to converge \\
 \hline
\end{tabular}}
\caption{Acronyms used to characterise simulation performance.}
\label{acro}
\end{table}

\subsection{Propagation of a single soliton}
If the initial condition is given by
\begin{equation}\label{ic1}
\psi_0(x) = \sqrt{\frac{2a}{q}} \exp\Bigg ( i\bigg (\frac{c(x-x_0)}{2}\bigg )\Bigg )\:
\text{sech}\big (\sqrt{a}(x-x_0)\big ),
\end{equation}
then the solution is given by the single soliton 
\begin{equation}\label{ss1}
\psi(x,t) = \sqrt{\frac{2a}{q}}\:
\exp\Bigg ( i\bigg (\frac{c(x-x_0)}{2}-\Big (\frac{{c^2}-4a}{4}\Big )t\bigg )\Bigg )\:
\text{sech}\big (\sqrt{a}(x-x_0-ct)\big ), 
\end{equation}
where $c$ is the speed of the soliton and $a$ is a real parameter which determines the
amplitude. The modulus $\rho=|\psi(x,t)|$ is characterised by a single soliton of amplitude 
$\sqrt{\frac{2a}{q}}$, which is located initially at $x=x_0$, and travels to the right with speed $c$. 

\begin{figure}[!ht]
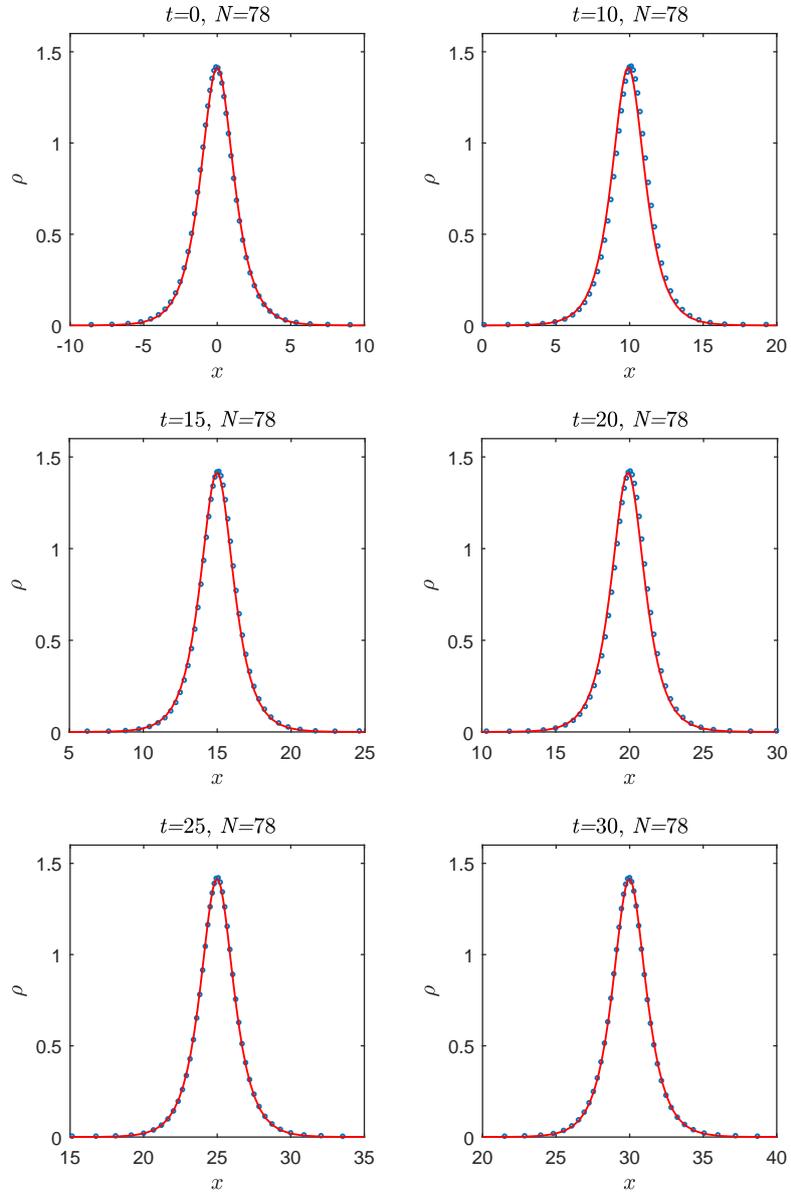

    \centering
    \scalebox{0.7}{\inc{w0-30}}
    \caption{Numerical solution (circles) and exact solution (red line) for single soliton.}
    \label{sch_sol1}
\end{figure}

\begin{figure}[!ht]
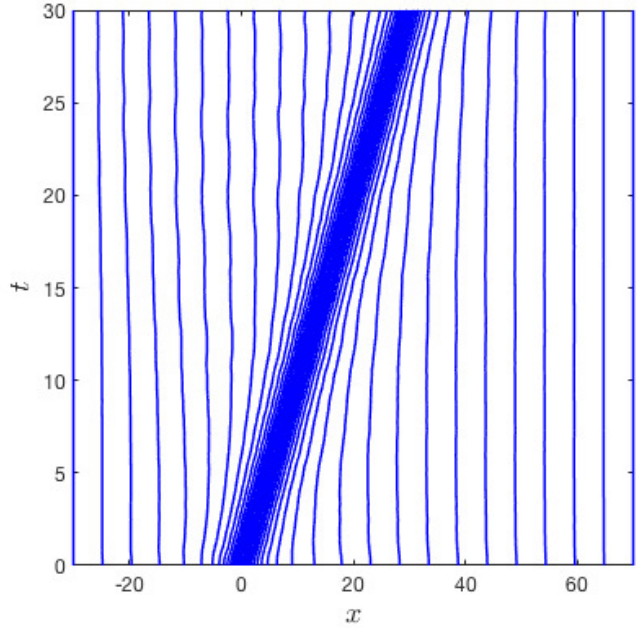

  \centering
  \scalebox{0.6}{\inc{traj}}
  \caption{Mesh trajectories for the propagation of a single soliton.}
  \label{s1traj}
\end{figure}

The computations shown here used the values $c=1$, $a=1$, $q=1$ and $x_0=0$.  For the spatial domain we set $x_l = -30$ and $x_r = 70$ and integrated up to time $T=30$. The simulations were performed using the following parameters and tolerances: $\ala=1.4$, $\bta=0.8$, \mbox{RTOL$=1.5\times10^{-2}$}, \mbox{ETOL$=5\times 10^{-3}$} and \mbox{$\tau=10^{-3}$}. 
Based on equidistribution of the initial condition, the algorithm determined that $N^0=78$ grid nodes were needed to resolve the initial condition to the required tolerance. Plots of the numerical and exact solutions at $t=0, 10, 15, 20, 25$ and $30$ are shown in Figure~\ref{sch_sol1}. We can see that the algorithm does an excellent job of resolving the travelling soliton. 
The grid point trajectories shown in Figure~\ref{s1traj} (showing alternate grid points) indicate that the $r$-adaptive component of the adaptive algorithm performs well tracking the movement of the  front throughout the time integration period.
The performance of the $hr$-algorithm is  shown in Table~\ref{tsch1}, where the acronyms used are described in Table~\ref{acro}. We can see that no $h$-refinement steps were needed throughout the integration process as the travelling wave was well resolved using the $r$-refinement component of the adaptive strategy. 

\begin{figure}[!ht]
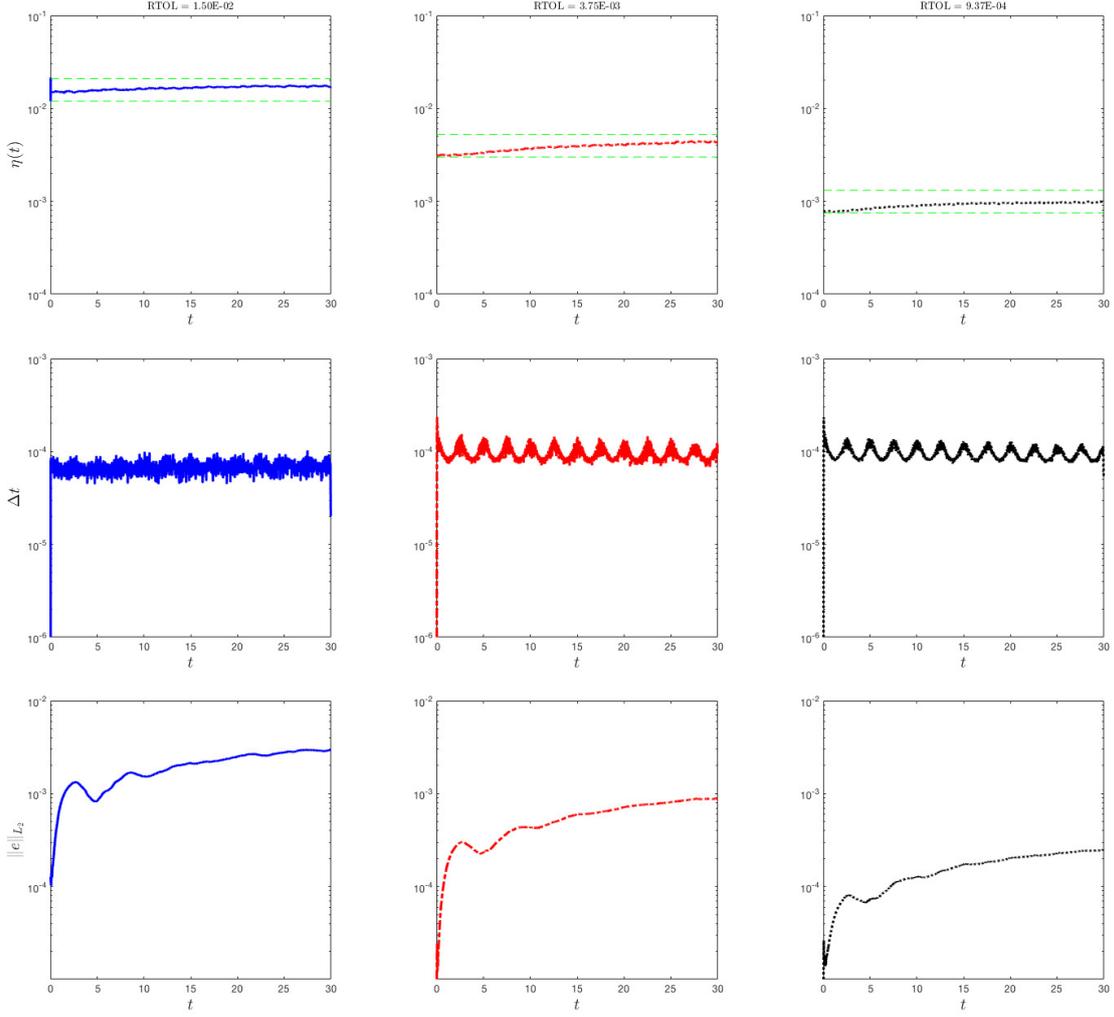

  \centering
  \scalebox{0.3}{\inc{eta2}} \hspace*{-1.2cm}
  \scalebox{0.3}{\inc{eta3}}  \hspace*{-1.2cm}
  \scalebox{0.3}{\inc{eta4}} \\
  \scalebox{0.3}{\inc{dt2}} \hspace*{-1.2cm}
  \scalebox{0.3}{\inc{dt3}} \hspace*{-1.2cm}
  \scalebox{0.3}{\inc{dt4}} \\
  \scalebox{0.3}{\inc{l22}} \hspace*{-1.2cm}
  \scalebox{0.3}{\inc{l23}} \hspace*{-1.2cm}
  \scalebox{0.3}{\inc{l24}}
 \caption{Evolution of $\eta$, $\Dta t$ and $||e||_{L_2}$ with time as RTOL is quartered using ETOL$=10^{-8}$.}
\label{RTOL}
\end{figure}

\begin{table}[!ht]
  \begin{center}
  {\renewcommand{\arraystretch}{1.5}
    \begin{tabular}{|cccccccc|}
      \hline
      NHR & NMAX & NMIN & NSTP & \hspace{.5cm}BS\hspace{.5cm} & \:CTF\: & \:ETF\: & \:JACS\: 
      \\ \hline 
      $0$ &$78$ &$78$ & $434$ &$2538$ &$0$ &$0$& $2538$ \\ \hline
    \end{tabular}}
  \end{center}
\caption{ Performance of the $hr$-algorithm for the propagation of a single soliton.}
\label{tsch1}
\end{table}
To compare our results with previous work, we computed the $L_2$ error at time $t^n$ as follows:
\begin{equation}\label{l2err}
||e^n||_{L_2} = \sqrt{\frac{1}{x_r-x_l}\sum_{i=1}^{N} \left(\frac{x_{i}^n - x_{i-1}^n}{2}
    \right) \left((e_i^n)^2 + (e_{i-1}^n)^2\right)   \:},
\end{equation}
where $e_i^n = \tilde{\rho}(x_i,t^n) - \rho(x_i, t^n)$ and $\tilde{\rho}$ is the exact solution. The first column of Figure~\ref{RTOL} shows the evolution of the spatial error indicator $\eta$, the time step history, and the $L_2$ norm of the error when a value ${\rm RTOL}=6\times 10^{-2}$ is used. One can see that the algorithm ensures that  
$\beta \:\rm{RTOL}\leq \eta(t)\leq \alpha \:\rm{RTOL}$, as expected. In addition, the time step is chosen to be reasonably large and is relatively constant throughout the simulation. One can also see that the error remains approximately constant throughout the simulation.

The second and third columns of Figure~\ref{RTOL} show how $\eta$, the time step, and the $L_2$ norm of the solution error at $t=T$, vary as RTOL is decreased by a factor of four.  We can see from Table \ref{tolprop} 
that $N^{0}$ is approximately doubled as ${\rm RTOL}$ is quartered and that this results in a reduction in 
the $L_{2}$-norm of the error and the rate of decrease is approaching a factor of four. This behaviour is consistent with the expectation that the overall algorithm is spatially second-order accurate. 

\begin{table}[!h]
\centering
{\renewcommand{\arraystretch}{1.5}
\begin{tabular}{|c|c|c|}
	\hline RTOL & $N^0$ & $||e||_{L_2}$ \\
	\hline $\phantom{00}1.5\times 10^{-2}$ & $67$ & $2.95\times 10^{-3}$ \\
	\hline $\phantom{0}3.75\times 10^{-3}$ & $175$ & $8.87\times 10^{-4}$ \\
	\hline $\phantom{0}9.37\times 10^{-4}$ & $352$ & $2.43\times 10^{-4}$ \\ \hline
\end{tabular}}
\caption{Tolerance proportionality of $hr$-algorithm for single soliton problem.}
\label{tolprop}
\end{table}

For this problem the exact values of the conserved quantities are $Q=4$ and \mbox{$E=-1/3$}. The conserved quantities $Q$ and $E$ are approximated at time $t^n$ by
\beq 
Q_h ^n  
= \sum_{i=1}^{N-1} \left(\frac{h_i^n + h_{i+1}^n}{2}\right) \left( (U_i^n)^2 + (V_i^n)^2
\right) \:,
\eeq
and
\begin{align*}
E_h ^n &= \sum_{i=1}^{N-1} h_i^{n}\Biggl[ \left(\frac{U_{i+1}^n -
	U_{i}^n}{h_i^n} \right)^2 +
\left(\frac{V_{i+1}^n-V_{i}^n}{h_i^n} \right)^2 
-\frac{q}{2}\left((U_i^n)^2+ (V_i^n)^2\right)^2 \Biggr].
\end{align*}
We define the mean values $\bar{Q}_h$ and $\bar{E}_h$ obtained from the numerical solution over the time integration period as
\beq 
\bar{Q}_h= \frac{1}{NT}\sum_{n=0}^{NT} Q_h^n \quad\quad {\rm and} \quad\quad \bar{E}_h= \frac{1}{NT}\sum_{n=0}^{NT} E_h^n,
\eeq
where $NT$ is the total number of time steps taken. Table~\ref{vsuni} shows the initial values $Q_{h}^{0}$ and $E_{h}^{0}$, and the 
mean values $\bar{Q}_{h}$ and $\bar{E}_{h}$ using the $hr$-adaptive algorithm. The equivalent quantities are also shown when two uniform meshes are used: one using the SDIRK2 time integration scheme and 
the other using the energy-conserving modified Crank-Nicolson scheme (MCN) \cite{del81}. We have also shown the 
error in $Q$ and $E$ when a finer uniform mesh is used. We can see that the initial approximations $Q_{h}^{0}$ and $E_{h}^{0}$ are much more accurate using the  
adaptive mesh compared to the schemes using a uniform mesh with the same number of mesh points; this is especially true for the approximation of $E$ which relies on the accurate approximation of solution  derivatives. We can see that there is no noticeable difference in the approximation to $Q$ and $E$ 
using the conservative scheme; the error being equal to its initial error as expected. We see that the value of $Q^{n}_{h}$ is almost constant throughout the computation using the $hr$-algorithm. The conservation of $E$ is less well maintained but $\bar{E}_h$ is 
far more accurately predicted with the $hr$-algorithm compared to a uniform mesh with the same number of mesh nodes.  
It requires at least $N=400$ nodes and a time step of size $\Dta t = 5\times 10^{-3}$ to achieve a comparable level of solution accuracy (in the $L_2$ norm) to an adaptive grid method, which utilises $N=78$ nodes and $\Dta t \approx 10^{-1}$. 
\begin{table}[!ht]
  \centering
  {\renewcommand{\arraystretch}{1.5}
\begin{tabular}{|c|c|c|c|c|c|c|}
	\hline
Approach &$N$ & $\Dta t$ & $|Q_h^0-Q|$ & $|\bar{Q}_h-Q|$ & $|E_h^0-E|$ & $|\bar{E}_h-E|$\\ \hline
$hr$-adaptivity & $78$ & $\approx 10^{-1}$ &  $7.9\times 10^{-3}$ & $7.0\times 10^{-3} $ & $1.4\times 10^{-2}$ & $2.8\times 10^{-2} $ \\ \hline
Uniform grid & $78$ & $5\times 10^{-3}$ & $4.5\times 10^{-2}$ &  $4.5 \times 10^{-2}$ & $3.2\times 10^{-1}$ & $3.2\times 10^{-1}$ \\ \hline
MCN & $78$ & $5\times 10^{-3}$ & $4.5\times 10^{-2}$ &  $4.5 \times 10^{-2}$ & $3.2\times 10^{-1}$ & $3.2\times 10^{-1}$ \\ \hline
Uniform grid & $400$ & $5\times 10^{-3}$ &  $ 5.3\times 10^{-15}$ & $ 4.6\times 10^{-7}$ & $2.1\times 10^{-2}$ & $2.1\times 10^{-2}$ \\ \hline
\end{tabular}} 
\caption{Comparison of conserved quantities using $hr$-adaptivity and a fixed uniform mesh.} 
\label{vsuni}
\end{table}

\begin{figure}[!ht]
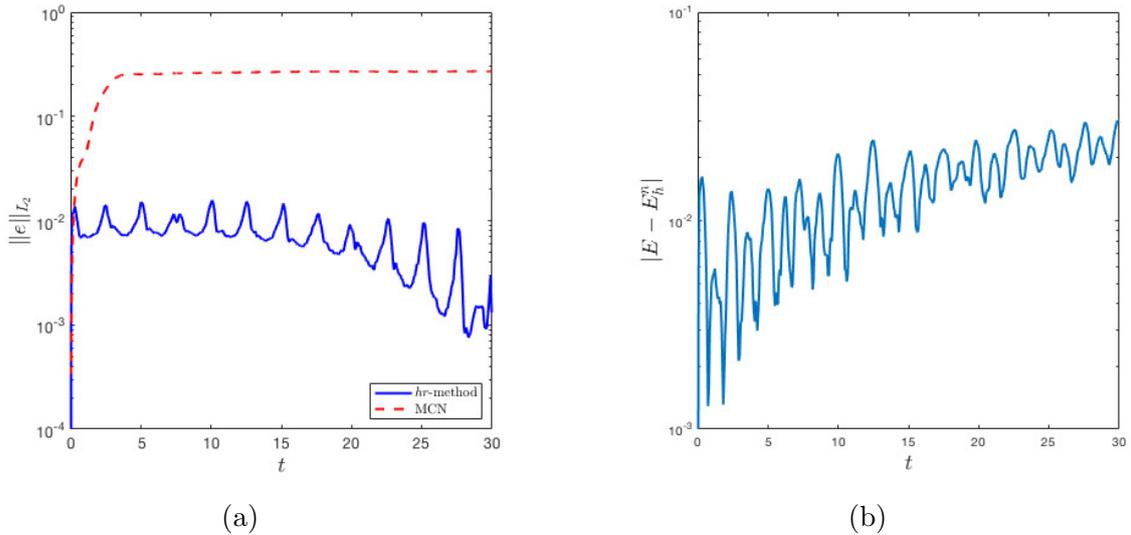

	\centering
	\begin{subfigure}{.5\textwidth}
		\scalebox{0.45}{\inc{l2vsMCN}}
		\caption{}
		\label{l2vsMCN}
	\end{subfigure}%
	\begin{subfigure}{.5\textwidth}
		\scalebox{0.45}{\inc{ErrEHT30}}
		\caption{}
		\label{c1Eerr}
	\end{subfigure}
	\caption{(a) Error in $L_{2}$-norm using the $hr$-adaptive and energy conserving MCN schemes. 
		(b) Error in the energy $|E - E_h^n|$ using the $hr$-adaptive method.}
	\label{c1Errcomp}
\end{figure}

In Figure \ref{l2vsMCN}, we plot the evolution of $L_{2}$-norm of the error. For comparison we have also included the error calculated using the MCN scheme on a uniform mesh with $N=78$. It is clear that the error using the moving mesh method is considerably smaller than the error obtained using the MCN scheme and this is despite the fact that 
the energy is not conserved with the moving mesh method as can be seen in Figure \ref{c1Eerr}. 
This example therefore suggests that it is important to accurately resolve spatial gradients and that it is 
not sufficient to use a numerical method which conserves an approximation to the energy if the spatial mesh is not dense enough. 

Table~\ref{terr} shows the $L_2$ error at $T=1$ for various values of $N$ using the moving mesh method when mesh enrichment and coarsening are switched off. These results are comparable with those obtained by previous authors \cite{rev:86, ssc:86, svs96}. The results obtained by Revilla \cite{rev:86} for a similar experiment are also displayed.   Whereas the $hr$-algorithm we use utilises interpolation only when the mesh is coarsened/enriched, the method used by Revilla obtains the numerical solution on the new grid via an interpolation based method. It is possible that frequent interpolation is a source of spatial error in the Revilla approach which could account for the discrepancy between the results. The approach in \cite{ssc:86} uses the arc-length monitor function based on first-derivatives of the numerical solution and interpolation to transfer the solution between meshes. It is clear that our method outperforms 
both of these moving mesh approaches and is far more accurate than the use of a uniform mesh.

\begin{table}[!ht]
  \centering
  {\renewcommand{\arraystretch}{1.5}
  \begin{tabular}{|c|c|c|c|c|}
    \hline
     & \multicolumn{4}{c|}{$||e||_{L_2}$} \\ \cline{2-5}
    \raisebox{1.5ex}[0pt]{$N$} & $r$-method & Revilla \;\cite{rev:86}\: & \cite{ssc:86} & Uniform\\ \hline
    $50$ & $1.8\times 10^{-3}$ & $1.8\times 10^{-2}$  & $6.8\times 10^{-2}$ & $9.4\times 10^{-2}$  \\ \hline
    $100$ & $4.6\times 10^{-4}$ & $3.3\times 10^{-3}$ & $1.6\times 10^{-2}$ & $6.1\times 10^{-2}$\\ \hline
    $200$ & $1.2\times 10^{-4}$ & $2.1\times 10^{-3}$& $3.9\times 10^{-3}$ & $1.4\times 10^{-2}$\\ \hline
  \end{tabular}}
\caption{$L_2$ error at $T=1$ using $r$-method for single soliton problem.} 
\label{terr}
\end{table}

\subsection{Interaction of two solitons}
We next consider the initial condition
\beq\label{ic2}
\psi(x,0) = \psi_1(x) + \psi_2(x), 
\eeq
where 
\bo
\psi_1(x) &=& \sqrt{\frac{2a_{1}}{q}} \exp\Bigg ( i\bigg (\frac{c_{1}(x-x_{01})}{2}\bigg )\Bigg )\:
\text{sech}\big (\sqrt{a_{1}}(x-x_{01})\big ), \\
\psi_2(x) &=& \sqrt{\frac{2a_{2}}{q}} \exp\Bigg ( i\bigg (\frac{c_{2}(x-x_{02})}{2}\bigg )\Bigg )\:
\text{sech}\big (\sqrt{a_{2}}(x-x_{02})\big ).
\eo
In the simulations below we have chosen the parameters  $a_1=0.2$, $c_1=1$, $x_{01}=0$ and $a_2=0.5$, $c_2=-0.2$, $x_{02}=25$ for $\psi_1$ and $\psi_2$, respectively with $q=1$. This initial condition consists of two solitons of different amplitudes, located initially at $x_{01}=0$ and $x_{02}=25$ which then move to the right and left with speeds $c_1$ and $c_2$, respectively. The solitons interact as if they were particle-like entities, i.e., they exhibit elastic collisions from which they emerge with the same shape \cite{molbook}.

\begin{figure}[!ht]
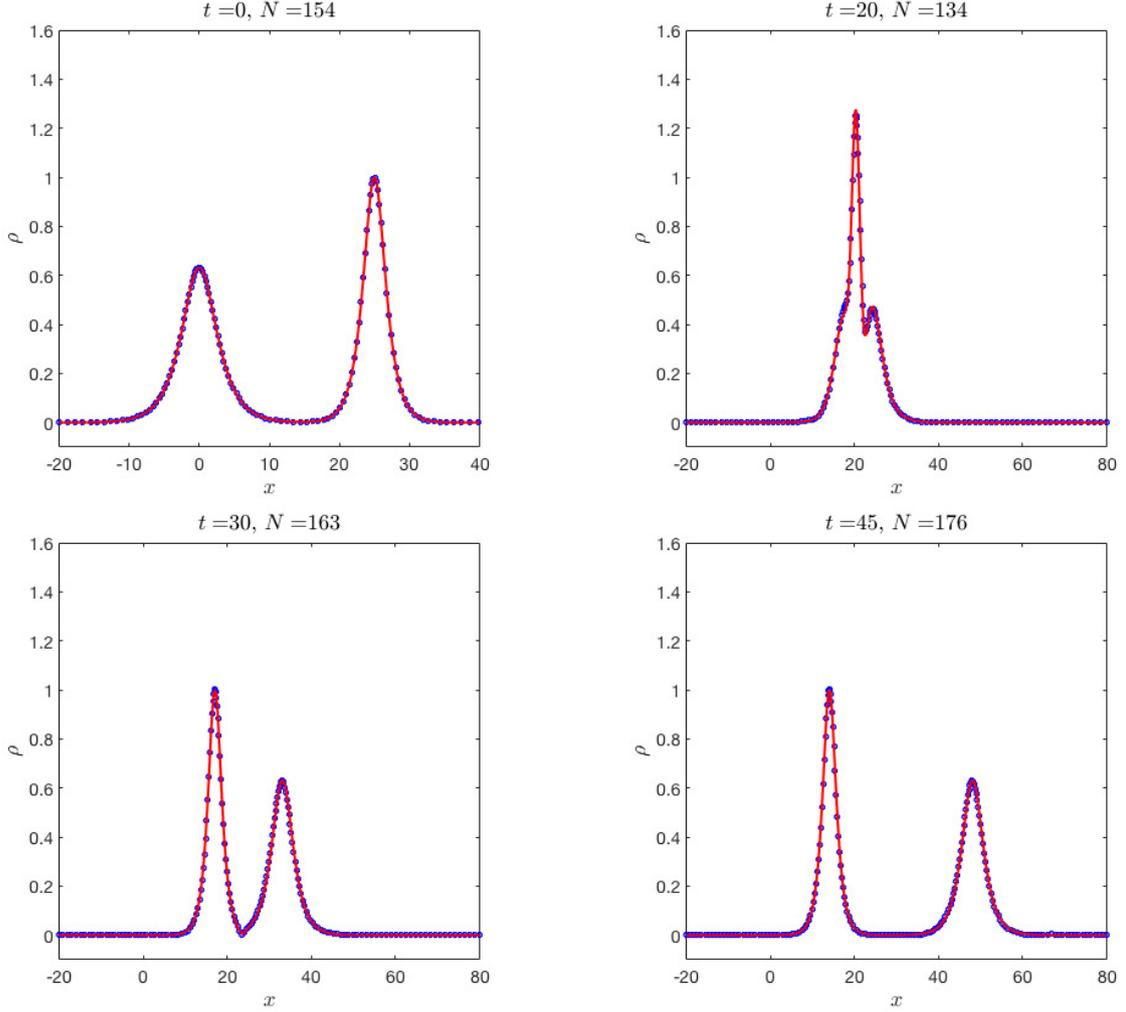

	\centering
	\begin{subfigure}{.5\textwidth}
		\centering
		\scalebox{0.45}{\inc{t0}}
	\end{subfigure}%
	\begin{subfigure}{.5\textwidth}
		\centering
		\scalebox{0.45}{\inc{t20}}
	\end{subfigure}
\\
	\centering
	\begin{subfigure}{.5\textwidth}
		\centering
		\scalebox{0.45}{\inc{t30}}
	\end{subfigure}%
	\begin{subfigure}{.5\textwidth}
		\centering
		\scalebox{0.45}{\inc{t45}}
	\end{subfigure}
	\caption{Numerical solution (circles) and reference solution (red line) for the interaction of two solitons using the adaptive $hr$-adaptive algorithm.} 
\label{sch_sol2}
\end{figure}

We integrated the problem up to the final time $T=45$ and set the artificial boundaries \mbox{$x_l=-20$} and $x_r = 80$. 
For this computation we used ETOL$=5\times 10^{-4}$, $\tau=10^{-2}$,  RTOL$=1\times 10^{-2}$, $\ala=1.2$ and $\bta=0.8$. The solutions at $t=0, 20, 30, 45$ are shown in Figure~\ref{sch_sol2}. The reference solution was obtained using a fixed uniform grid with $2000$ nodes. To plotting accuracy, the $hr$-algorithm has done an excellent job of capturing the interacting solitons. Details of the performance of the $hr$-algorithm are given in Table~\ref{tsch2}. 

\begin{table}[!ht]
  \begin{center}
  {\renewcommand{\arraystretch}{1.5}
    \begin{tabular}{|cccccccc|}
      \hline
      NHR & NMAX & NMIN & NSTP & \:JACS\: & \hspace{.5cm}BS\hspace{.5cm} & \:ETF\: & \:CTF\: 
      \\ \hline 
      $7$ &$197$ &$134$ & $909$ &$4524$ &$4524$ &$0$& $0$ \\ \hline
    \end{tabular}}
  \end{center}
\caption{Performance of $hr$-algorithm for the interaction of two solitons.}
\label{tsch2}
\end{table}

\begin{figure}[!h]
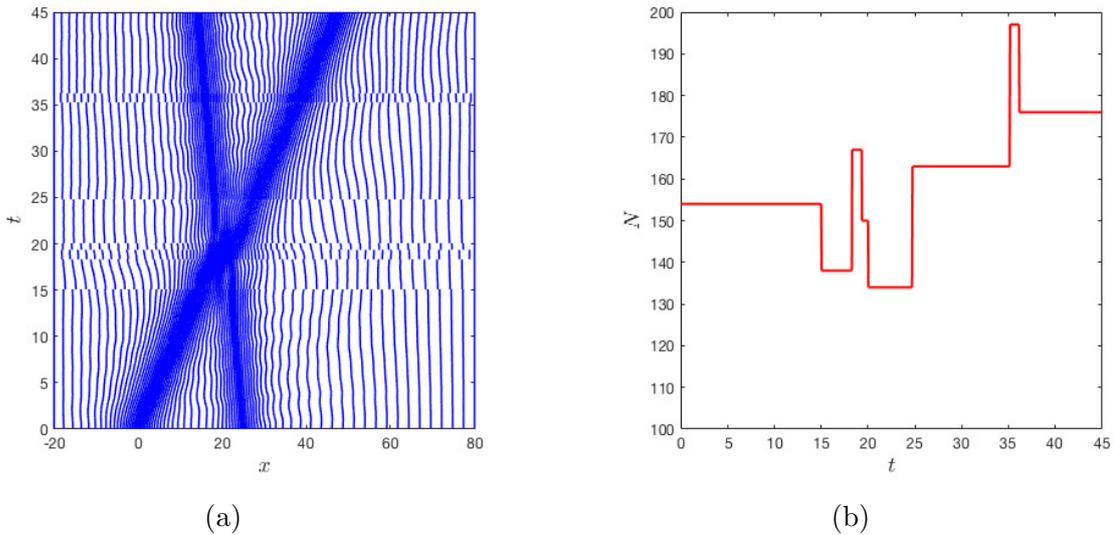

	\centering
	\begin{subfigure}{.5\textwidth}
		\scalebox{0.45}{\inc{traj2}}
		\caption{}
		\label{c2traj}
	\end{subfigure}%
	\begin{subfigure}{.5\textwidth}
		\scalebox{0.45}{\inc{N}}
		\caption{}
		\label{c2N}
	\end{subfigure}
	\caption{(a) Mesh trajectories and (b) evolution of $N$ for the two soliton problem.}
	\label{2soltraj}
\end{figure}

In Figure~\ref{c2traj} we can see from the grid point trajectories that between $t=15$ and $t=25$ the two solitons meet and cross. It is within this time interval that refinements/derefinements of the grid are performed as can be seen from the grid point trajectories and evolution of $N$ in Figure~\ref{c2N}. 
Note that the behaviour of $N$ in Figure~\ref{2soltraj} is what one would anticipate. The amplitudes of the solitons at $t=0$ and $t=45$ are about the same and so we would expect roughly the same number of nodes to be used. We can see that the algorithm utilised a very similar value for $N$ at these times. Furthermore, more nodes were introduced just as the problem got more challenging and nodes were removed when they were not needed.

\begin{figure}[!ht]
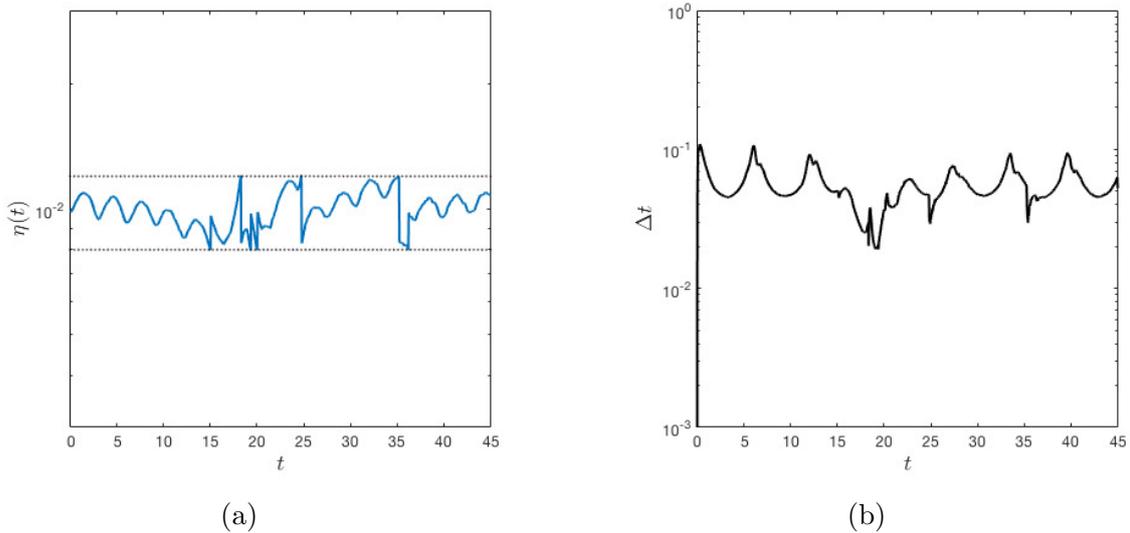

	\centering
	\begin{subfigure}{.5\textwidth}
		\scalebox{0.45}{\inc{eta}}
		\caption{}
		\label{s2eta}
	\end{subfigure}%
	\begin{subfigure}{.5\textwidth}
		\scalebox{0.45}{\inc{dt}}
		\caption{}
		\label{s2dt}
	\end{subfigure}
	\caption{(a) Evolution of $\eta$ and (b) time step history for the two soliton interaction problem.}
	\label{2soleta}
\end{figure}

\begin{figure}[!ht]
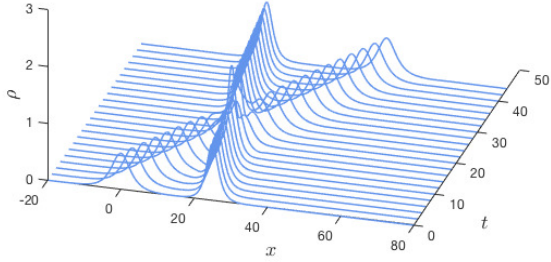
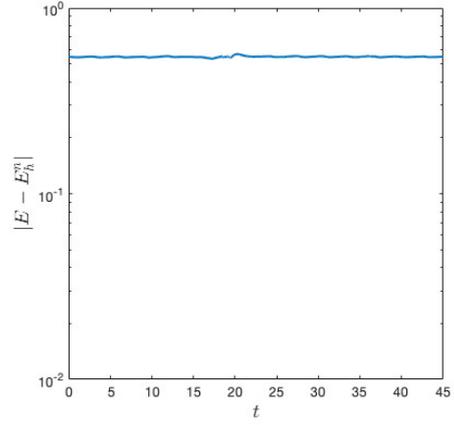

	\centering
	\begin{subfigure}{.5\textwidth}
		\scalebox{0.4}{\inc{c2profile}}
		\caption{Solution using $hr$-method}
		\label{c2pro}
	\end{subfigure}%
	\begin{subfigure}{.5\textwidth}
		\scalebox{0.40}{\inc{E_err}}
		\caption{Energy error.}
		\label{c2Eerr}
	\end{subfigure}
	\caption{(a) Solution profiles using the $hr$-method and (b) the evolution of the error in energy for two interacting solitons where the exact energy $E=-2/3$}
	\label{c2profiles}
\end{figure}

It is clear from Figure~\ref{s2eta} that the objective of keeping $\eta$ within the given bounds has been successfully achieved. Also, we see from Figure~\ref{s2dt} that the algorithm has chosen the time step to be fairly constant throughout the simulation - it reduces appropriately when the solitons interact and recovers afterwards.  A global plot of the solution profiles is given in Figure \ref{c2profiles} as well as a plot of the error in the energy. We can see that there is a slight change in the energy as the solitons interact but this error does not increase dramatically through the simulation and is not far off its initial value.

\begin{figure}[!h]
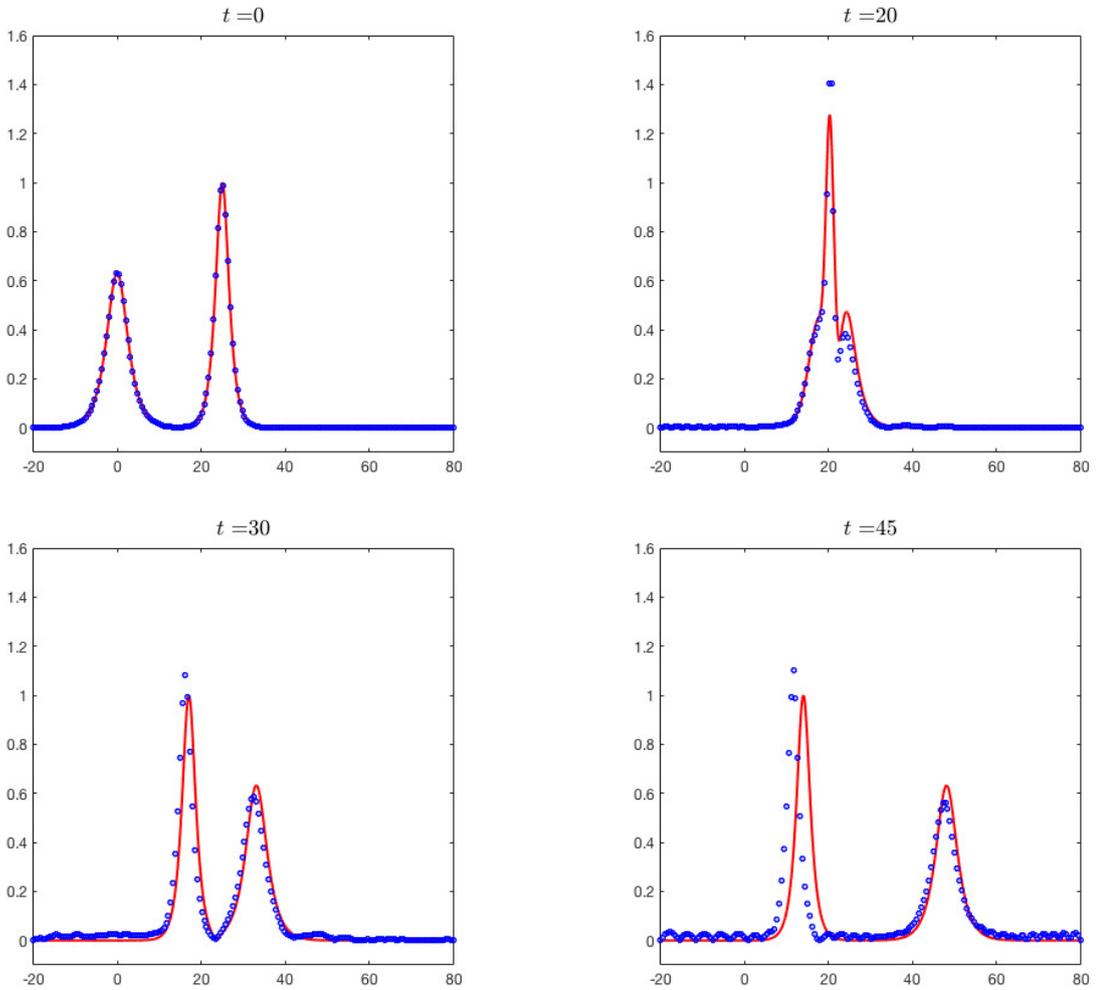

	\centering
	\begin{subfigure}{.5\textwidth}
		\scalebox{0.45}{\inc{mcn177t1}}
	\end{subfigure}%
	\begin{subfigure}{.5\textwidth}
		\scalebox{0.45}{\inc{mcn177t2}}
	\end{subfigure}
\\
	\centering
	\begin{subfigure}{.5\textwidth}
		\scalebox{0.45}{\inc{mcn177t3}}
	\end{subfigure}%
	\begin{subfigure}{.5\textwidth}
		\scalebox{0.45}{\inc{mcn177t4}}
	\end{subfigure}
	\caption{Numerical solution (circles) and reference solution (line) for interacting solitons using MCN scheme and a uniform mesh with $N=177$.} 
	\label{sch_sol2MCN}
\end{figure}

For comparison in Figure \ref{sch_sol2MCN} we show the computed solution obtained using the 
MCN scheme with a time step $\Delta t=10^{-3}$  and a uniform mesh using the average number of points used with the adaptive $hr$ scheme. For this example this resulted in a value of $N=177$. We can see clearly that the computed solutions are poor in comparison with the $hr$-adaptive method. The inaccuracy again is due to the lack of spatial resolution of the rapidly changing behaviour of the interacting solitions and this is not compensated for by the discrete conservation of energy. 

\subsection{Bound state of three solitons}
The bound state of multiple solitons is a class of problems with initial condition
\beq\label{u03}
\psi_0(x) = \sech (x), 
\eeq
and $q=2N_{s}^2$ in \eqref{sch}, where $N_{s}$ is a positive integer. For this case, $\psi(x,t)$ is a bound state of $N_{s}$ solitons \cite{m81}. 

As in \cite{hmm85,rev:86, svs96} we study the bound state of three solitons, i.e., $q=18$, and use artificial boundaries located at $x_l=-20$ and $x_r=20$. The solution of this problem is periodic in time  with the period being approximately $T=0.8$. The solution develops extremely large spatial and temporal gradients thus posing a stringent test to any numerical scheme. 

For the purpose of comparison with previous work, we integrated over five periods and took $T=4$. For this simulation, we  used the values RTOL$=10^{-3}$, $\ala=3$, $\bta=0.4$ and ETOL$=5\times 10^{-3}$. For this problem we set  $\varphi_{u,v}=10^{-3}$ throughout the computation. Note that the choice \eqref{eq:mu} and \eqref{eq:mv} places approximately half of the available nodes outside the region where we have high solution activity. For this problem, this leads to insufficient resolution in the soliton region. The value of $\varphi_{u,v}=10^{-3}$ is much smaller than that predicted by the algorithm thus placing more nodes in the centre of the domain. 

\begin{figure}[!ht]
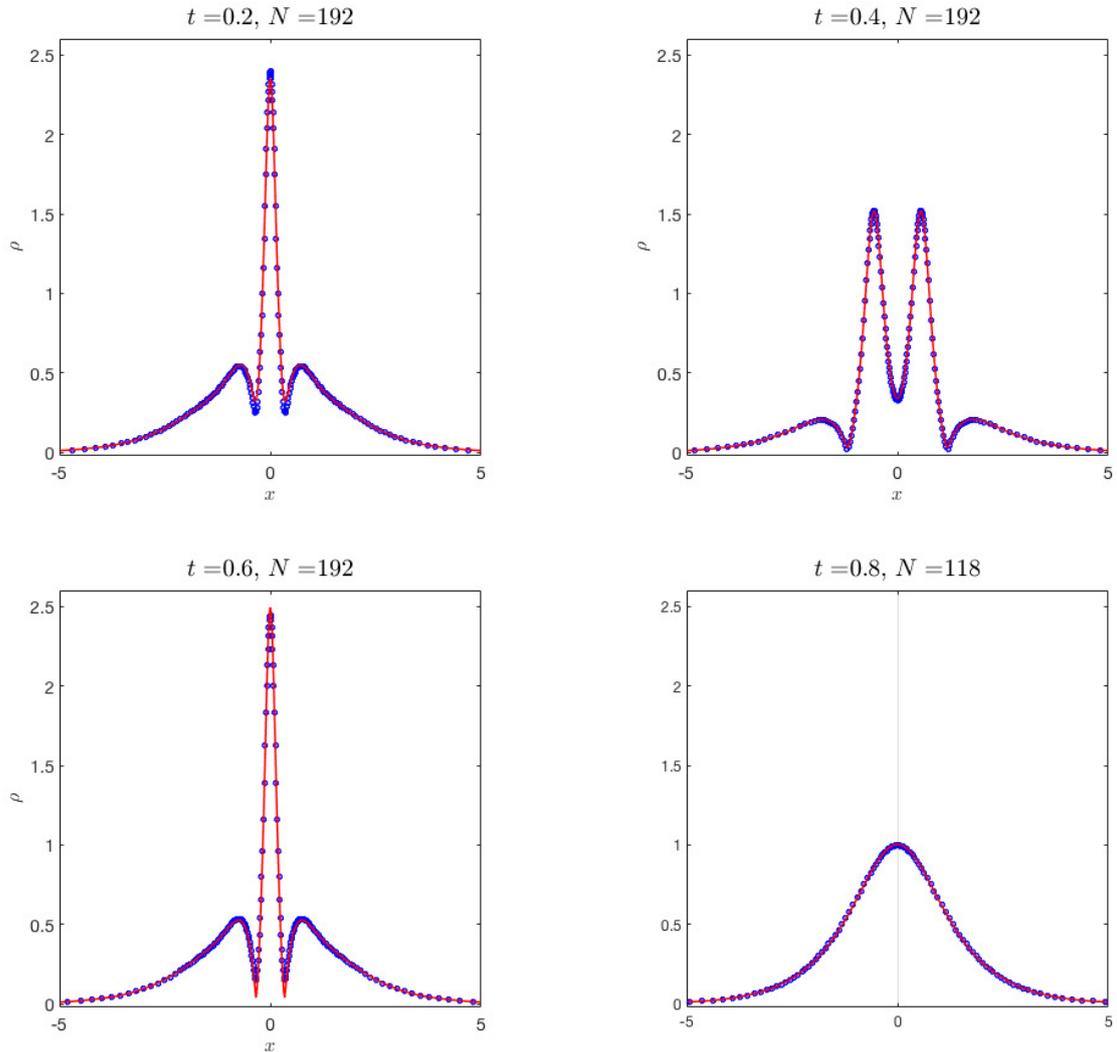

	\centering
	\begin{subfigure}{.5\textwidth}
		\scalebox{0.45}{\inc{p4t1}}
		\label{s31}
	\end{subfigure}%
	\begin{subfigure}{.5\textwidth}
		\scalebox{0.45}{\inc{p4t2}}
		\label{s32}
	\end{subfigure} 
\\
	\begin{subfigure}{.5\textwidth}
		\scalebox{0.45}{\inc{p4t3}}
		\label{s33}
	\end{subfigure}%
\begin{subfigure}{.5\textwidth}
	\scalebox{0.45}{\inc{p4t4}}
	\label{s34}
\end{subfigure} 
	\caption{Numerical solution (circles) and reference solution (line) for the bound state of three solitons using the adaptive $hr$-algorithm.} 
	\label{sch_sol3}
\end{figure}

\begin{figure}[!h]
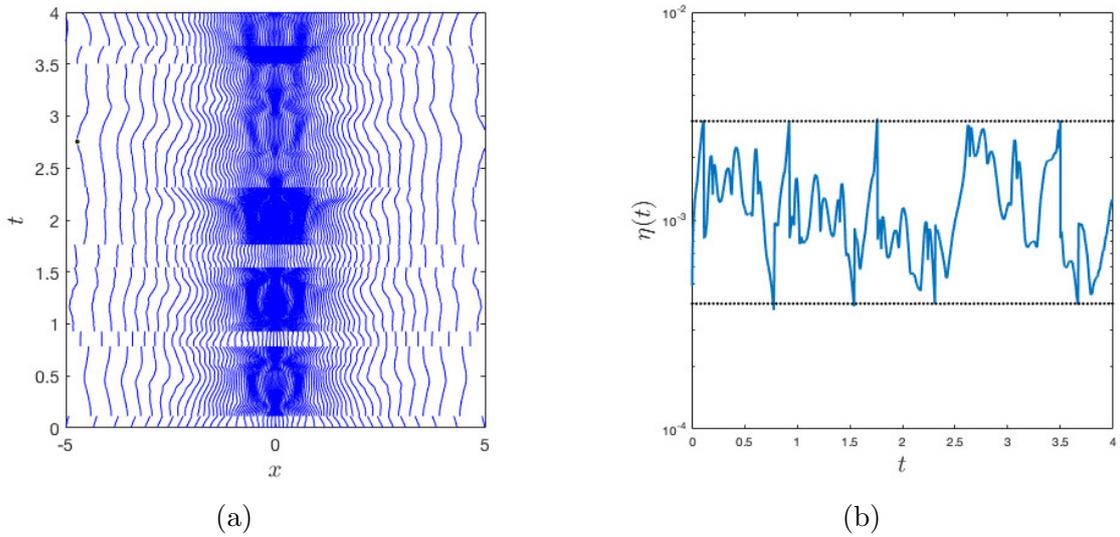

	\centering
	\begin{subfigure}{.5\textwidth}
		\scalebox{0.45}{\inc{p4trajz}}
		\caption{}
		\label{s3traj}
	\end{subfigure}%
	\begin{subfigure}{.5\textwidth}
		\scalebox{0.45}{\inc{p4eta}}
		\caption{}
		\label{s3eta}
	\end{subfigure}
	\caption{(a) Mesh trajectories for the region $[-5,5]\times [0,1]$ and (b) evolution of $\eta$ for bound state of three solitons.}
	\label{3sol}
\end{figure}

The solutions at the times $t=0.2, 0.4, 0.6$ and $0.8$ are shown in Figure~\ref{sch_sol3}. The reference solution was obtained using a fixed uniform grid with $N=2000$. The mesh trajectories and the evolution of $\eta$ are displayed in Figure~\ref{3sol}. The trajectories only cover the region $-5 \leq x \leq 5$ where nodes experience significantly more activity. As for the previous cases, the algorithm tracks accurately the solitons as they appear and disappear during the integration. One can also see that the monitor captures regions of high curvature as expected. The error indicator $\eta$ is also kept within the required bounds. Computational statistics are shown in Table~\ref{tsch3}.

\begin{table}[!ht]
	\begin{center}
		\begin{tabular}{|cccccccc|}
			\hline
			NHR & NMAX & NMIN & NSTP & \:JACS\: & \hspace{.5cm}BS\hspace{.5cm} & \:ETF\: & \:CTF\: 
			\\ \hline 
			$8$ &$332$ &$100$ & $856$ &$5157$ &$5157$ &$7$& $0$ \\ \hline
		\end{tabular}
	\end{center}
	\caption{Computational statistics for bound state of three solitons over 5 periods.}
	\label{tsch3}
\end{table}

\clearpage

\begin{figure}[!h]
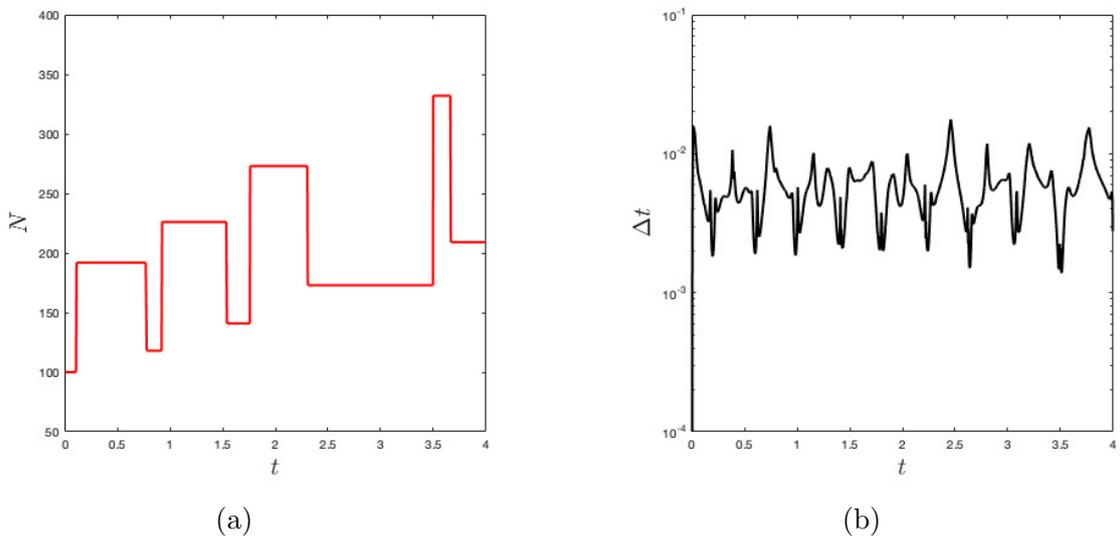

	\centering
	\begin{subfigure}{.5\textwidth}
		\scalebox{0.45}{\inc{p4N}}
		\caption{}
	\end{subfigure}%
	\begin{subfigure}{.5\textwidth}
		\scalebox{0.45}{\inc{p4dt}}
		\caption{}
	\end{subfigure}
	\caption{(a) Evolution of $N$ and (b) time step history for bound state of three solitons.}
	\label{3ndt}
\end{figure}

\begin{figure}[!ht]
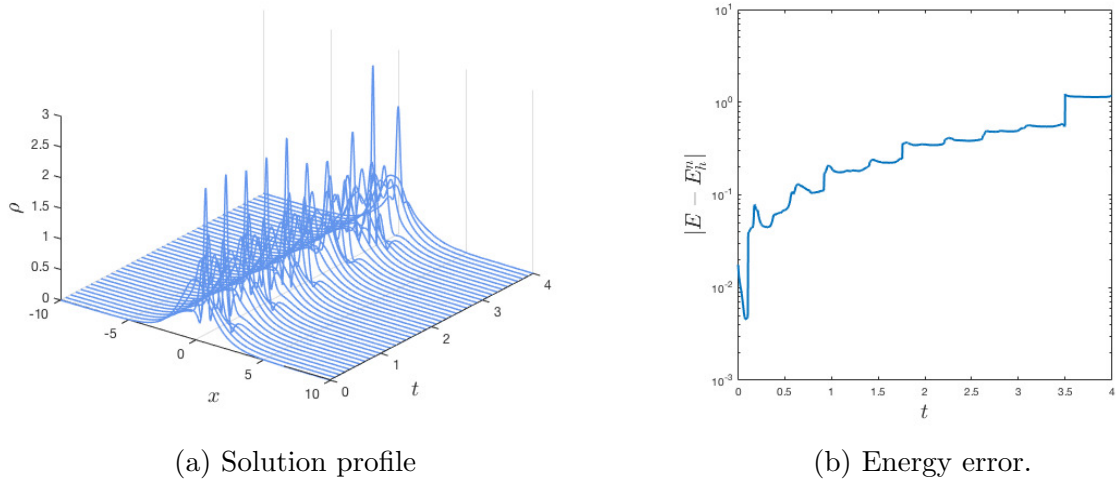

	\centering
	\begin{subfigure}{.5\textwidth}
		\scalebox{0.4}{\inc{P4profile}}
		\caption{Solution profile}
	\end{subfigure}%
	\begin{subfigure}{.5\textwidth}
		\scalebox{0.4}{\inc{E_rrT4}}
		\caption{Energy error. }
	\end{subfigure}
\caption{(a) Solution profile over 5 periods using the $hr$ method and (b) the error in the energy where the exact energy $E=-34/3$.}
\label{c3profile}
\end{figure}

\begin{figure}[!ht]
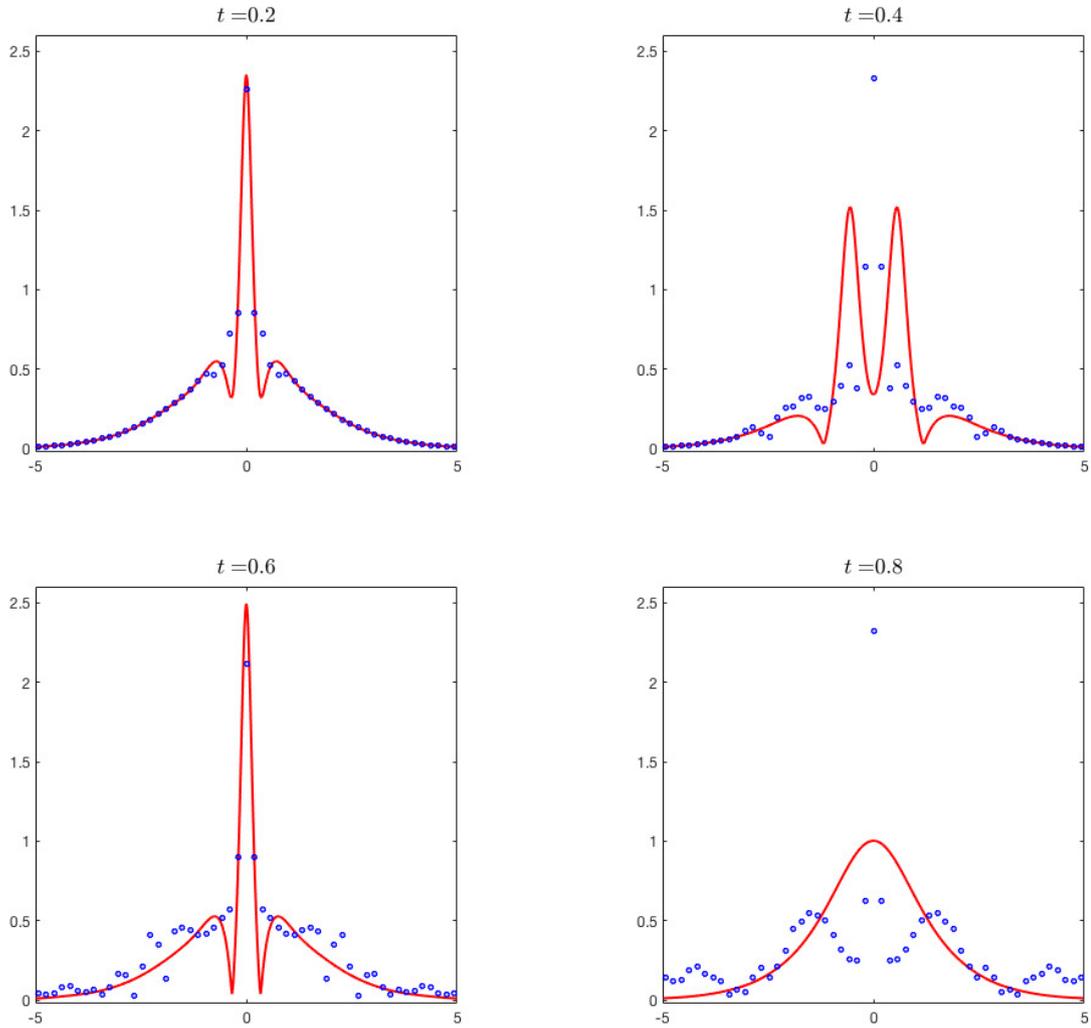

	\centering
	\begin{subfigure}{.5\textwidth}
		\scalebox{0.45}{\inc{mcn210t1}}
		\label{s31m}
	\end{subfigure}%
	\begin{subfigure}{.5\textwidth}
		\scalebox{0.45}{\inc{mcn210t2}}
		\label{s32m}
	\end{subfigure} 
\\
	\begin{subfigure}{.5\textwidth}
		\scalebox{0.45}{\inc{mcn210t3}}
		\label{s33m}
	\end{subfigure}%
	\begin{subfigure}{.5\textwidth}
		\scalebox{0.45}{\inc{mcn210t4}}
		\label{s34m}
	\end{subfigure} 
	\caption{Numerical solution (circles) and reference solution (line) for the bound state of three solitons using MCN scheme and a uniform mesh with $N=210$.} 
	\label{sch_sol3_mcn}
\end{figure}

The evolution of $N$ and the time step history are shown in Figure~\ref{3ndt}. As expected, the algorithm adds nodes when they are needed (in this case when there are two solitons) and removes them when they are not needed. Furthermore, the time step is kept relatively large over the integration, another benefit of using an adaptive approach. 
A global view of the solution trajectories is shown in Figure~\ref{c3profile} 
as well as the evolution of the energy error. We can see that the solution is well approximated over the first few periods but that eventually the solution accuracy deteriorates. The energy error also increases more rapidly for this much more demanding test case as can be seen also in Figure~\ref{c3profile}. However, the computed solutions are far more accurate to those obtained using the MCN scheme on a uniform mesh as seen in Figure~\ref{sch_sol3_mcn}, where again we have used the average number of points used with the $hr$ method and a time step $\Delta t=10^{-4}$. As with the two previous examples, it's clear that a lack of spatial resolution results in poor solution accuracy even though the scheme preserves the inital approximation to the energy. 
\section{Conclusions}
\label{sec:conc}
We have developed an $hr$-adaptive algorithm which combines adaptive mesh movement and mesh enrichment for the cubic nonlinear Sch\"{o}dinger equation. Numerical experiments demonstrate that the method works well for problems with large spatial and temporal variations in comparison with other 
moving mesh methods as well as energy-conserving uniform grid methods with the equivalent number of grid nodes. 

There is scope to improve the algorithm further and we hope to do this by considering more closely the effect and the choice of MMPDEs and monitor functions. It would also be useful to develop conservative versions of the moving mesh method which could potentially improve solution accuracy over longer time integration periods. For simplicity, in this work we have used linear interpolation to transfer solutions between meshes when $h$-adaptation is required. In future we will also consider the use of $L_{2}$ projection 
for transferring solutions between meshes. While we have focussed in this paper on the application of method to the cubic NLSE, it also has applicability to other non-linear dispersive wave equations. Another area of development would be to investigate the merits of using an $hr$-adaptive  scheme in a multidimensional context.

\bibliographystyle{plain} 
\bibliography{PaperRefs}
\end{document}